\documentclass[12pt]{article} 
\usepackage{amssymb,amsmath}
\usepackage{graphicx}
\usepackage{hyperref}
\usepackage{ifthen}

\makeatletter
\def\@seccntformat#1{\csname the#1\endcsname.\ } 
\def\@biblabel#1{#1.} 
\makeatother

\newif\ifNoRemark
\def\addtheorem#1#2#3#4{
\ifthenelse{\equal{#2}{}}{}%
{\ifthenelse{\expandafter\isundefined\csname the#2\endcsname}{\newcounter{#2}}{}}
\newenvironment{#1}[1][\global\NoRemarktrue]
{\par\addvspace{2mm plus 0.5mm minus 0.2mm}\noindent 
\ifthenelse{\equal{#2}{}}{}{\refstepcounter{#2}}%
{\bf #3%
\ifthenelse{\equal{#2}{}}{}{{~\csname the#2\endcsname}}%
 \vphantom{##1}\ifNoRemark.\ \else\ (##1).\fi}\begingroup #4}%
{\endgroup\par\addvspace{1mm plus 0.5mm minus 0.2mm}\global\NoRemarkfalse}
\expandafter\newcommand\csname b#1\endcsname{\begin{#1}}
\expandafter\newcommand\csname e#1\endcsname{\end{#1}}
}

\addtheorem{theorem}{thrm}{Theorem}{\sl}
\addtheorem{lemma}{lmm}{Lemma}{\sl}
\addtheorem{corollary}{corollary}{Corollary}{\sl}
\addtheorem{example}{example}{Example}{}
\addtheorem{remark}{remark}{Remark}{\sl}

\newenvironment{proof}[1][\hspace{-1.0ex}]%
{\par\addvspace{1mm}{\sc Proof\hspace{1.0ex}{#1}.} }%
{\quad$\blacktriangle$\par\addvspace{1mm}}
\def\qed{}

\title{On weight distributions of perfect colorings and completely regular codes%
\thanks{This is an author version of the paper published in the Designs Codes and Cryptography,
DOI~\href{http://dx.doi.org/10.1007/s10623-010-9479-4}{10.1007/s10623-010-9479-4},
\copyright 2010 Springer}
\thanks{The results of the paper were partially presented at the
Sixth International Workshop on
Optimal Codes and Related Topics
in June 2009, Varna, Bulgaria.}
}
\author{\href{http://arXiv.org/a/krotov_d_1}{Denis S. Krotov}%
\thanks{Sobolev Institute of Mathematics, pr. Akad. Koptyuga 4, Novosibirsk 630090, Russia.
\texttt{krotov@math.nsc.ru}}
}
\date{}

\def\trn{^{\rm\scriptscriptstyle T}}
\begin{document}
\maketitle
\begin{abstract}
A vertex coloring of a graph is called ``perfect'' if for any two colors $a$ and $b$,
the number of the color-$b$ neighbors of a color-$a$ vertex $x$ does not depend on
the choice of $x$, that is, depends only on $a$ and $b$ (the corresponding partition
of the vertex set is known as ``equitable'').
A set of vertices is called ``completely regular'' if the coloring according
to the distance from this set is perfect. By the ``weight distribution''
of some coloring with respect to some set we mean the information
about the number of vertices of every color at every distance from the set.

We study the weight distribution of a perfect coloring
(equitable partition) of a graph with respect to a completely regular set
(in particular, with respect to a vertex if the graph is distance-regular).
We show how to compute this distribution by the knowledge of
the color composition over the set.
For some partial cases of completely regular sets, we derive explicit formulas
of weight distributions. Since any (other) completely regular set itself generates
a perfect coloring, this gives universal formulas for calculating the weight distribution
of any completely regular set from its parameters.
In the case of Hamming graphs,
we prove a very simple formula for the weight enumerator
of an arbitrary perfect coloring.
\def\keywords#1{\par Keywords:~#1 \par}
\end{abstract}

\section{Introduction}

A remarkable property of perfect codes in Hamming spaces is that the weight distribution
of the code with respect to some vertex depends only on the distance between this vertex
and the code and the parameters of the code
(i.e., the code minimal distance, the dimension of the space, and the size of the alphabet) \cite{Lloyd,ShSl}.
One of the quite general generalizations of the perfect codes that inherit this property is
the perfect colorings (also known as equitable partitions or regular partitions), or their
real-valued generalizations, which are also considered in this paper and named perfect structures.
In particular, the coloring of the vertices according
to the distance from a perfect code is a perfect coloring.
If any (not necessarily perfect) code generates a perfect coloring in such a way,
then it is called a completely regular code.

In this paper we consider a matrix way to calculate weight distributions of perfect structures
with respect to completely regular codes (not only sole vertices).
We derive general formulas that include as parameters the parameters of the perfect structure
and the completely regular code.
As pointed out in Section~\ref{s:dr},
the matrix formulas for calculating the weight distribution of
a perfect coloring (equitable partition) with respect to a vertex
of a distance-regular graph was known before \cite{Martin:PhD},
but it seems that they have never been published in a journal or considered in the general form
presented in the current paper.

The paper is organized as follows.
The first part, Sections~\ref{s:PC,CR}--~\ref{s:we},
contain general facts and formulas for calculating weight distributions.
In Section~\ref{s:PC,CR}, we define the perfect colorings,
completely regular sets and perfect structures, and consider several examples.
In Section~\ref{s:distr}, we define the distribution of one perfect structure with respect to another
and discuss two simple algebraic facts (Theorems~\ref{th:distrib} and~\ref{th:main}) that have
important corollaries for perfect structures.
In Section~\ref{s:dr}, we show how to compute the weight distribution of a perfect structure with respect
to a point in distance-regular graphs.
In Section~\ref{s:we}, for the case of Hamming graphs $H_q^n$,
we prove a simple formula for the weight enumerator
of an arbitrary perfect structure.

The second part, Sections~\ref{s:somesets} and~\ref{s:local},
is devoted to more special areas, where we consider distributions
with respect to some sets and local distribution of perfect structures, mainly in Hamming graphs.
This part, while being not so interesting from the theoretical point of view, provides
potentially useful tools for studying and characterizing different kinds of perfect structures.
In Section~\ref{s:somesets}, we consider weight distributions with respect to some special completely regular sets,
mainly in Hamming spaces.
In Section~\ref{s:local}, we consider so-called local distributions:
given a coloring $f$ of the cartesian product $G'\times G''$
of two graphs, a distribution of $f$ in some instance of, say, $G''$ is called a local distribution.
It turns out that all such local distributions form a perfect structure over $G'$; this allows one to derive some
relations between them.
As an example, we calculate local distributions for perfect codes in $H_{10}^{11}$,
in the case of their existence.

\section{Perfect colorings and completely regular sets}\label{s:PC,CR}

Let $G=(V=\{0,\ldots,N-1\},E)$ be a graph;
let $ f$ be a function (``coloring'') on $V$ that possesses exactly ${\nu}$ different values
$e_0$, \ldots, $e_{{\nu}-1}$ (``colors'').
The function $ f$ is called a \emph{perfect coloring} with parameter matrix $S=(S_{ij})_{i,j=0}^{\nu-1}$, or \emph{$S$-perfect coloring},
if
for any $i,j$ from $0$ to ${\nu}-1$ any vertex of color $e_i$ has exactly
$S_{ij}$ neighbors of color $e_j$.
(The corresponding partition of $V$ into ${\nu}$ parts is known
as an \emph{equitable partition}. In another terminology, see, for example, \cite{CDS},
$f$ is called an \emph{$S$-feasible  coloration} and $S$ is called a \emph{front divisor} of $G$.)

In what follows, we assume that $e_i$ is the tuple with $1$ in the $i$th position and $0$s in the others
(the length of the tuple may vary depending on the context; in the considered case it is ${\nu}$).
Denote by $A$ the adjacency matrix of $G$.
Then it is easy to see \cite[Lemma~5.2.1]{Godsil93}
that $ f$ is an $S$-perfect coloring
if and only if
\begin{equation}\label{eq:Af=fS}
Af = fS
\end{equation}
where the function $f$ is represented by its value array;
that is, the $i$th row of the $|V|\times {\nu}$ matrix $f$
is $f(i)$.
If the equation (\ref{eq:Af=fS}) holds for some matrices
$A$, $S$, and $f$ (of size $N\times N$, ${\nu} \times {\nu}$, and $N \times {\nu}$ respectively) over
$I\hspace{-0,9ex}R$,
then we will say that $f$ is an \emph{$S$-perfect structure}
(or a perfect structure with parameters $S$)
 over $A$ \cite{Avg:lect}.

So, in this context, the concept of perfect structure is a continuous generalization of the concept of
perfect coloring.
Conversely, a perfect coloring (equitable partition)
is equivalent to a perfect structure over the graph (i.e., over its adjacency matrix)
with the rows from $e_0$, \ldots, $e_{{\nu}-1}$.

Suppose that $f$ satisfies (\ref{eq:Af=fS}) with a three-diagonal parameter matrix $S$.
In this case, the corresponding perfect coloring (if any)
has the following property:
the colors $e_i$, $e_j$ of any two neighbor vertices satisfy $|i-j|\leq 1$.
The support of the $e_0$ of such coloring is known as a \emph{completely regular set}, or
\emph{completely regular code}
with \emph{covering radius} ${\nu}-1$. In other words, a set $C$ of vertices of a graph $G=(V,E)$ is
a \emph{completely regular set} if and only if its \emph{distance coloring}
(i.e., the function $ f(x) = e_{d_G(x,C)}$ where
$d_G(\cdot,\cdot)$ is the natural distance in the graph)
is perfect.

In the rest of this section, we consider several examples of perfect structures.
For examples of perfect $2$-colorings
of binary $n$-cubes, see \cite{FDF:PerfCol,FDF:12cube.en};
of Johnson graphs, see \cite{AvgMog:J63J73};
of halved $n$-cubes, see \cite{Kro:24}. For examples of completely regular codes
in $q$-ary $n$-cubes, see \cite{RifZin07}; in Johnson graphs, see \cite{Mar:CRD}.
\begin{example}
Very partial, but also a very important
case of perfect structures is the case of ${\nu}=1$;
then $f$ is just an eigenvector of $A$
(in the graph case, an eigenfunction of the graph)
with the eigenvalue equal to the only element of $S$.
\end{example}

\begin{example}
A graph is \emph{distance-regular} if the distance coloring with respect to any vertex is perfect with parameters
that do not depend on the choice of the vertex. An equivalent definition of distance-regular graphs
is given in Section~\ref{s:dr}.
\end{example}

\begin{example}\label{ex:1perf}
A subset $P$ of the vertex set $V$ of a regular graph $G=(V,E)$ is
known to be a \emph{$1$-perfect code}, if its distance coloring is perfect with
the parameter matrix {\small$\left(\begin{array}{cc}0&d\\1&d{-}1\end{array}\right)$},
where $d$ is the degree of the graph.
$1$-perfect codes in $n$-cubes (see Example~\ref{ex:cube} below) are actively studied;
the best-investigated case is binary, see, for example, \cite{Hed:2008:survey,Sol:2008:survey},
but even in that case the problem of full characterization of such codes is
far from being solved.
%

\end{example}
\begin{example}\label{ex:n-2}
A code $C$ in the binary $n$-cube
has cardinality $2^n/(n+3)$ and
minimal distance $3$ between codewords
(i.e., the parameters of doubly-shortened
$1$-perfect code)
if and only if it is the support
of the first color of a perfect coloring
with parameters
{\small$
\left(
\begin{array}{cccc}
0 & 1 & n{-}1 & 0 \\
1 & 0 & n{-}1 & 0 \\
1 & 1 & n{-}4 & 2 \\
0 & 0 & n{-}1 & 1
\end{array}
\right)
$}
\cite{Kro:2m-3}.
This gives an interesting example of non-completely-regular codes whose
code parameters guarantee that the code is a color of
some perfect coloring.
\end{example}
\begin{example}\label{ex:n-3}
An even more interesting example is the case of the codes with
parameters of triply-shortened $1$-perfect codes.
Such codes in general cannot be represented as a color of a perfect coloring
because the weight distribution of such a code depends on the choice of the code and
the choice of the initial codeword. Nevertheless, all the codes with these parameters
can be characterized in terms of perfect structures:
the distance coloring of such a code $C$ together with the distance coloring of its
antipode $C+11...1$ form a perfect structure with $6\times 6$ parameter matrix,
see \cite{Kro:2010ACCT:n-3} for details.
\end{example}
\begin{example}
A code $C$ in the binary $n$-cube
is called \emph{Preparata-like}
if it has cardinality $2^{n-1}/(n+1)^2$, $n=2^{2m}-1$, and
minimal distance $5$ between the codewords.
Equivalently, $C$ is the support of the first color of a perfect coloring with parameters
{\small$
\left(
\begin{array}{cccc}
0 & n & 0 & 0 \\
1 & 0 & n{-}1 & 0 \\
0 & 2 & n{-}3 & 1 \\
0 & 0 & n & 0
\end{array}
\right)
$}
(the parameters of the perfect coloring follow from the results of \cite{SZZ:1971:UPC}).
It is notable that unifying the first and the fourth colors of such a coloring
results in a $1$-perfect code, see Example~\ref{ex:1perf}.
\end{example}
\begin{example}\label{ex:sts}
A  $t$-$(n,k,\lambda)$ \emph{design} is an $n$-element set $V$ together with a set of $k$-element subsets of $V$ (called blocks) with the property that every $t$-element subset of $V$ is contained in exactly $\lambda$ blocks.
A  $(k{-}1)$-$(n,k,\lambda)$ design is equivalent to a perfect coloring of the Johnson graph $J(n,k)$
(see Example~\ref{ex:J} below) with parameters
{\small$
\left(
\begin{array}{cc}
k(\lambda{-}1) & k(n{-}k{-}\lambda{+}1)\\
k\lambda & k(n{-}k{-}\lambda)
\end{array}
\right),
$}
where the blocks are the vertices of the first color.
The blocks of any $(k{-}2)$-$(n,k,1)$ design also form a completely regular code in $J(n,k)$
\cite[Corol.~3.8]{Mar:CRD},
but with covering radius $2$.
\end{example}
Examples~\ref{ex:1perf}--\ref{ex:sts} show that some classes of codes or designs
with specified parameters
can be alternatively defined in terms of perfect structures (not necessarily
perfect colorings).
It seems quite important and useful to represent different
known classes of objects in terms of
perfect structures, and Examples~\ref{ex:n-2} and~\ref{ex:n-3}
show that this is possible even in the case of non-completely-regular sets.
It would be very interesting to find another example of such kind.
Natural candidates for consideration from this point of view are
so-called uniformly packed codes \cite{BZZ:1974:UPC}, or
different subclasses of such codes. The definitions
do not imply direct connections between uniformly packed codes
(in the sense of \cite{BZZ:1974:UPC})
and perfect structures, and finding this connection even in partial
cases seems to be a nontrivial problem.
One of intriguing examples of such a problem is to represent
the so-called Goethals-like codes (see \cite{ZinHel04} for recent results and bibliography)
 in terms of perfect structures.

\section{Distributions}\label{s:distr}

Assume that we have two perfect structures $f$ and $g$ over $A$
with parameters $S$ and $R$, respectively.
Then we say that $g\trn f$ is the \emph{distribution of $f$ with respect to $g$}.
For perfect colorings, it has the following sense: the element in the $i$th column
and $j$th row of $g\trn f$ equals the number of the vertices $v$ such that $g(v)=e_i$ and $f(v)=e_j$
(to avoid misunderstanding, we note that the number of elements in $e_i$
in the first equation is the number of colors of the coloring $g$,
while the number of elements in $e_j$ in the last equation is the number of colors of $f$;
in general, these numbers can be different, even if $i=j$).
In the case when $g$ is the distance coloring with respect to some (completely regular) set $C$,
we will also say that $g\trn f$ is the \emph{weight distribution of $f$ with respect to $C$}
(if $C=\{c\}$, with respect to $c$).
In other words, the weight distribution of $f$ with respect to $C$ is the tuple
$(h_0, h_1, h_2, ..., h_{\rho(C)} )\trn$ where $h_w$ is the sum of $f$
over all the vertices at the distance $w$ from $C$
and $\rho(C)=\displaystyle\max_{x \in V(G)}{d_G(x,C)}$ is the covering radius of $C$.

The two following theorems are elementary from an algebraic point of view;
nevertheless,
they are very significant for the perfect structures.
\begin{theorem}\label{th:distrib}
Let $f$ and $g$ be $S$- and $R$- perfect structures over $A$ and
$A^T$ respectively ($f$, $g$, $A$, $S$, and $R$ are $N\times {\nu}$,
$N \times {\mu}$, $N \times N$, ${\nu} \times {\nu}$, and
${\mu}\times {\mu}$ matrices). Then $g\trn f$ is a perfect structure
over $R\trn$ with parameters $S$. Briefly,
$$(Af=fS)\, \& \, (A\trn
g=gR) \Rightarrow  (R\trn g\trn f = g\trn f S).$$
\end{theorem}
\begin{proof} $R\trn g\trn f = (g R)\trn f = (A\trn g)\trn f = g\trn A f = g\trn f S$. \qed\end{proof}
Note that if $A$ is the adjacency matrix of some graph, then $A=A\trn$.

\begin{theorem}\label{th:main}
  If the matrix $B=\{B_{i,j}\}_{i,j=0}^{n-1}$ satisfies
  $B_{i,i+1}\neq 0$, $B_{i,j}= 0$ for any $i=0,...,n-2$, $j>i+1$, then any
  $S$-perfect structure $h$ over $B$ is uniquely defined by its first row $h_0$.
  Moreover, the rows $h_0, \ldots, h_{n-1}$ of $h$ satisfy the recursive relation
  \begin{equation}\label{eq:rec}
   h_{i+1} = (h_{i} S - \sum_{j=0}^{i} B_{i,j} h_j)/B_{i,i+1}, \qquad  i=0,\ldots,n-2,
   \end{equation}
  and, by induction,
\begin{equation}\label{eq:poly}
  h_i = h_0 \Pi_i^{(B)}(S)
\end{equation}
  where $\Pi_i^{(B)}(x)$ is a degree-$i$ polynomial in $x$.
\end{theorem}
\begin{proof}
From $Bh = hS$ we have $\sum_{j=0}^n B_{i,j} h_j = h_i S$.
Applying the hypothesis on $B_{i,j}$, we get (\ref{eq:rec}).
\qed\end{proof}

\begin{remark}
  In the important subcase of three-diagonal matrix $B$, the recursive relation 
  has the following form, where $B_{0,-1}$ and $h_{-1}$ are assumed to be zero:
\begin{equation}\label{eq:req3}
h_{i+1} = (h_{i} S - B_{i,i} h_{i} - B_{i,i-1} h_{i-1})/B_{i,i+1}, \qquad  i=0,\ldots,n-2.
\end{equation}
\end{remark}

So, given a completely regular set $C$, we also have a way to
reconstruct the weight distribution with respect to $C$ of any other
perfect structure (perfect coloring) $f$ over the same graph by
knowledge of only the first component of the distribution (the sum
of the function $f$ over the set $C$). To do this, we should apply
Theorem~\ref{th:main} with $B=R\trn$, where the three-diagonal
matrix $R$ is the parameter matrix of the distance coloring with respect
to $C$. The uniqueness of such reconstruction was known before
\cite{Avg:PhD}, but known formulas cover only partial cases of $f$,
for example, the weight distribution (with respect to a vertex) of
$1$-perfect binary codes can be found in \cite{Lloyd,ShSl}.

\section{Weight distributions in a distance-regular graph}\label{s:dr}

Let $G=(V,E)$ be a graph and let, for every $w$ from $0$ to $\mathrm{diam}(G)$
(the diameter of $G$), the matrix
$A^{(G)}_w=A_w=(a^w_{ij})_{i,j\in V}$
be the distance-$w$ matrix of $G$ (i.e., $a^w_{ij}=1$ if the graph distance between $i$ and $j$ is $w$,
and $a^w_{ij}=0$ otherwise); put $A=A_1$. The graph $G$ is called \emph{distance-regular} if,
for every $w$, the matrix $A_w$ equals $\Pi_w(A)$
for some  polynomial $\Pi_w$ of degree $w$.
The polynomials $\Pi_0$, $\Pi_1$, \ldots, $\Pi_{\mathrm{diam}(G)}$ are called \emph{$P$-polynomials} of $G$.

Now, suppose that $f$ is a perfect structure over $G$ (i.e., over $A$) with some parameters $S$.
By the definition, we have
\begin{equation}\label{eq:AffS}
A f = f S.
\end{equation}
From (\ref{eq:AffS}) we easily derive $A^w f = A^{w-1} f S = \ldots = f S^w$ for any degree $w$ and,
consequently, $P(A) f = f P(S)$ for any polynomial $P$.
In particular,
\begin{equation}\label{eq:AS}
A_w f = f \Pi_w(S).
\end{equation}
We now observe that the $i$th row of $A_w f$ is the sum of the vector-function $f$ over
all the vertices at distance $w$ from the $i$th vertex.
So, (\ref{eq:AS}) means the following:

\begin{theorem}\label{th:distr}
Assume we have an $S$-perfect structure $f$ over a distance-regular graph
with $P$-polynomials $\Pi_w$.
If the value of $f$ at a vertex $x$ is $f_0$, then
the tuple
$$(f_0=f_0 \Pi_0(S), f_0 \Pi_1(S), \ldots , f_0 \Pi_{\mathrm{diam}(G)}(S))\trn$$
is the weight distribution of $f$ with respect to $x$.
\end{theorem}
For a perfect coloring, the statement of the theorem means the following:
if the color of the vertex $x$ is $e_j$,
then the color composition of the vertices at distance
$w$ from $x$ is calculated as $e_j \Pi_w(S)$.
This fact (in terms of equitable partitions) was known before
\cite[Sect.\,2.2.2,\,2.1.5]{Martin:PhD}, and, from the algebraic point
of view, the generalization to perfect structures is not essential
and can hardly be considered as a new result. Nevertheless,
as we will see in Sections~\ref{s:somesets} and~\ref{s:local},
this generalization allows us to apply a common approach in studying
distributions with respect to several kind of subsets, not only sole vertices.

%
%

\begin{example}\label{ex:cube}
  Let $G=H^n_q$ be the \emph{$q$-ary $n$-cube},
  whose vertex set is the set of all
 $n$-words over the alphabet $\{0,\ldots,q-1\}$,
 two vertices being adjacent if and only if they differ in exactly one position.
 Then
  \begin{equation}\label{eq:prekrawtchuk}
  \Pi_w(\cdot) = P_w(P_1^{-1}(\cdot))
  \end{equation}
  where
  \begin{equation}\label{eq:krawtchuk}
   P_w(x) = P_w(x;n,q) = \sum_{j=0}^w (-1)^j (q-1)^{w-j} \left(x \atop j\right) \left(n-x \atop w-j\right)
  \end{equation}
 is the \emph{Krawtchouk polynomial}; $P_1(x) = (q-1)n-qx$.
A connected component of the distance-$2$ graph of the binary $n$-cube is a distance-regular graph
with $\Pi_w(\cdot) = P_{2w}(P_2^{-1}(\cdot))$ (recall $P_2(x) = \left(n \atop 2 \right) - 2nx +2x^2$),
known as the \emph{halved $n$-cube}.
\end{example}


\begin{example}\label{ex:J}
 Let $G=J(n,k)$ be the \emph{Johnson graph}, whose vertex set is the set of all binary
 $n$-tuples with exactly
 $k$ ones, two vertices being adjacent if and only if they differ in exactly two positions.
 Then
 $\Pi_w(\cdot) = E_w(E_1^{-1}(\cdot))$ where
   $$E_w(x) = E_w(x;n,k) = \sum_{j=0}^w (-1)^j \left(x \atop j\right) \left(k-x \atop w-j\right) \left(n-k-x \atop w-j\right)$$
 is the \emph{Eberlein polynomial} \cite{Delsarte:1973}.
\end{example}

\section{Weight enumerators in Hamming spaces}\label{s:we}

Assume that $(h_0, h_1, h_2, ..., h_n)\trn$ is the weight distribution
with respect to some fixed point
of a perfect structure $f$ over the $q$-ary $n$-cube.
By the \emph{weight enumerator} of $f$ we will mean the
vector-valued polynomial
$$
W_f(z)=h_0+zh_1+z^2h_2+\dots+z^nh_n
$$
in a real-valued variable $z$.
\begin{theorem}\label{th:we}
Let $f$ be an $S$-perfect structure over the $q$-ary $n$-cube $H_q^n$;
let $h_0$ be the value of $f$ in some fixed point.
Then the weight enumerator $W_f(z)$ of $S$ with respect to this point satisfies
\begin{eqnarray}\label{eq:we}
W_f(z)&=&h_0 Z_{S;n,q}(z) , \mbox{ where }\\ \label{eq:we0}
Z_{S;n,q}(z)&=&(1-z)^{((q-1)nI-S)/q}(1+(q-1)z)^{(nI+S)/q}.
\end{eqnarray}
\end{theorem}
\begin{proof}
1. We first consider the case when the rank of $f$ coincides with the size of $S$. It is known and easy to check that the Krawtchouk polynomials (\ref{eq:krawtchuk}) satisfy
\begin{equation}\label{eq:gp}
(1-z)^{x}(1+(q-1)z)^{n-x}=\sum_{w=0}^n P_w(x)z^w
\end{equation}
for every integer $x$ from $0$ to $n$. Taking into account
(\ref{eq:AS}) with the accompanying observation and (\ref{eq:prekrawtchuk}),
we have to prove that (\ref{eq:gp}) is true for $x$ equal to the matrix
$P_1^{-1}(S) = ((q-1)nI-S)/q$. To prove this, it is sufficient to show
that this matrix is diagonalizable and its
eigenvalues are integers from $1$ to $n$.
Equivalently, $S$ has a complete set of eigenvectors with eigenvalues
from $\{-n, -n{+}q, -n{+}2q, \ldots, (q{-}1)n\}$.
But this is true for the adjacency matrix $A$ of $H_q^n$ (see \cite[Theorem~9.2.1]{Brouwer}).
It is easy to see from $Af=fS$
that if $v$ is an eigenvector of $S$,
then $fv$ is an eigenvector of $A$ with the same eigenvalue;
so, the restrictions on the eigenvalues of $S$ are proved.
Moreover, if $u$ is a generalized eigenvector of $S$ and
$Su = \lambda u + v$, then $A(fu)= \lambda (fu)+(fv)$, i.e., $fu$
is a generalized eigenvector of $A$, which is impossible because $A$ is symmetric.
So, there are no generalized eigenvectors of $S$, and hence,
$S$ is diagonalizable.

2. Now, consider an arbitrary case. Let $\nu'$ be the rank of the matrix $f$. Then there are $N\times \nu'$ matrix $f'$, $\nu\times \nu'$ matrix $t'$ and $\nu'\times \nu$ matrix $t$ such that $ft'=f'$ and $f't=f$. From $Af=fS$ we derive
$Af'=f'S'$ with $S'=tSt'$; that is, $f'$ is an $S'$-perfect structure.
Since the rank of $f'$ coincides with the size of $S'$,
we can apply p.1 
to get
$$
W_{f'}(z)=h'_0(1-z)^{((q-1)nI-S')/q}(1+(q-1)z)^{(nI+S')/q}.
$$
where $h'_0=h_0 t'$ is the value of $f'$ in the initial point.
Since $\phi(S')=t\phi(S)t'$ for any analytical function $\phi$, we also have (\ref{eq:we}).
\qed\end{proof}

\begin{remark}
  If $f$ is a perfect coloring with $\nu$ colors, then the rank of the matrix
  $f$ is $\nu$; as a corollary, the eigenvalues of the parameter
  matrix are eigenvalues of the graph (we do not need p.2 of Theorem~\ref{th:we} in this case).
  Nevertheless, the last is not true for some perfect structures
  derived from perfect colorings in Sections~\ref{s:somesets}
  and~\ref{s:local}.
\end{remark}

\section{Distributions with respect to some sets}\label{s:somesets}

In this section, we will derive formulas for the weight distributions of perfect structures
with respect to some special completely regular sets, which have large covering radius and small
($1$ or $2$) code distance. As we will see, for the considered cases, the situation is reduced to
calculating the weight distributions
with respect to a vertex in some smaller distance-regular graph.

\subsection{A lattice}

The set $R$ discussed in this subsection plays some role in the theory
of perfect colorings of the $n$-cube.
It occurs in constructions of perfect colorings \cite{FDF:PerfCol,FDF:12cube.en};
it necessarily occurs in any linear distance-$2$ completely regular binary code \cite{Vas:cr2};
in particular, in shortened $1$-perfect binary code, and a (nonshortened) variation of this set (case $m=2$), known as a linear $i$-component, is widely used for the construction of $1$-perfect binary codes, see, for example, \cite{Sol:2008:survey}.
We will derive a rather simple formula for the weight distribution of a perfect coloring with respect to $R$.
The name ``lattice'' written in the title of the section comes from some attempts to draw
such a set in a figure.

Let us consider the $q$-ary $mk$-cube $H^{mk}_q$ and the
function $ \tilde g : V(H^{mk}_q) \to V(H^{k}_q)$ defined as
\begin{equation}\label{eq:g}
 \tilde g(x_1,\ldots,x_m)=x_1+\ldots + x_m \bmod q, \qquad x_i\in H^{k}_q .
\end{equation}
The set $R$ is defined as the set of zeroes of $\tilde g$.

\begin{lemma}\label{l:mA}
  $\tilde g$ is a perfect coloring of $H^{mk}_q$ with the matrix $m A^{(H^{k}_q)}$ where
$A^{(H^{k}_q)}$ is the adjacency matrix of $H^{k}_q$.
\end{lemma}
\begin{proof}
Consider the neighborhood of a vertex $x=(x_1,\ldots,x_m)$ of color
$y = x_1+\ldots + x_m \bmod q$. It is sufficient to show that it contains
 $m$ vertices of every color $z$ adjacent to $y$.
Indeed, all these vertices are obtained from $x$ by adding $z-y \bmod q$
to one of its components $x_1,\ldots,x_m$. \qed
\end{proof}

So, after representing the values of $\tilde g$
by the corresponding tuples $e_i \in \{0,1\}^{V(H^{k}_q)}$
$$ g(x) = e_{\tilde g(x)}, $$
we have the equation
$$ A^{(H^{mk}_q)} g =  g m A^{(H^{k}_q)}.  $$
By Theorem~\ref{th:distrib},
for any other perfect structure $f$ over $H^{mk}_q$ with parameter matrix $S$,
we have
$ (m A^{(H^{k}_q)}) (g\trn f) = (g\trn f) S $
or, equivalently,
\begin{equation}
\label{eq:gfmS}
 A^{(H^{k}_q)} (g\trn f) = (g\trn f) \left(\frac 1m S\right).
\end{equation}
That is, $(g\trn f)$ is a perfect structure over $H^{k}_q$ with parameters $\frac 1m S$.
Taking into account the following simple fact, we see that our problem is reduced to the calculation of
the weight distribution of this new perfect structure with respect to the zero vertex.
\begin{lemma}\label{l:dist}
The distance from a vertex $x$ to $R$ coincides
with the distance from $\tilde g(x)$ to the zero.
\end{lemma}
\begin{proof}
Clearly, modifying $x$ in one position, we vanish at most one element
of $\tilde g(x)$. On the other hand, as follows from Lemma~\ref{l:mA},
at least one (indeed, any) element can be vanished in such a way. So, the statement
of the lemma is proved by induction on the number of nonzero elements in $\tilde g(x)$.
\qed
\end{proof}

So, we can use the results of the previous sections to calculate
the weight distribution of $f$ with respect to $R$.

\begin{theorem}\label{th:grid}
Let $f$ be an $S$-perfect structure over the $q$-ary $mk$-cube $H^{mk}_q$.
Let $f_0$ be the sum of $f$ over the set $R$ of zeroes of $\tilde g$ {\rm(\ref{eq:g})}.
Then the weight distribution of $f$ with respect to $R$ is
$$ (f_0=f_0\Pi_0(\textstyle\frac 1m S),\, f_0 \Pi_1(\frac 1m S),\, \ldots ,\,
f_0 \Pi_{k}(\frac 1m S))\trn,$$
where $\Pi_i(\cdot)= P_i(P_1^{-1}(\cdot))$,
$P_i(\cdot) = P_i(\cdot;k,q)$,
see {\rm(\ref{eq:prekrawtchuk})} and {\rm(\ref{eq:krawtchuk})}.
The corresponding weight enumerator
equals $f_0 Z_{\frac 1m S;k,q}(z)$ where $Z_{\ldots}$ is defined in {\rm(\ref{eq:we0})}.
\end{theorem}
\begin{proof}
As follows from Lemma~\ref{l:dist},
the weight distribution of $f$ with respect to $R$
coincides with the weight distribution of $(g\trn f)$
with respect to the zero vertex.
Since, by (\ref{eq:gfmS}), $(g\trn f)$ is a $\frac 1m S$-perfect structure over $H^{k}_q$,
Theorem~\ref{th:distr} gives the required formula for the weight distribution
(for the explicit formula of the $P$-polynomials see Example~\ref{ex:cube}).
The formula for the weight enumerator comes from Theorem~\ref{th:we}.
\qed
\end{proof}

\subsection{The cartesian product}
Here, we consider distributions with respect to an instance of one of the multipliers
in the cartesian product of two graphs. Given two graphs $G'=(V',E')$, $G''=(V'',E'')$,
their cartesian product $G' \times G'' = (V,G)$ is defined as follows: the vertex set
is the set $V'\times V'' = \{(v',v'') : v'\in V', v'' \in V'' \}$; two vertices
$u=(u',u'')$ and $v=(v',v'')$ are adjacent (i.e. $\{u,v\}\in E$) if and only if
either $u'=v'$ and $\{u'',v''\} \in E''$ or $\{u',v'\} \in E'$ and $u''=v''$.

Let us consider the projection of $G' \times G''$ into $G''$:
\begin{equation}\label{eq:h}
\tilde h(x',x'')=x'', \qquad x'\in V', \qquad x''\in V''.
\end{equation}
Let us fix some vertex $o$ (say, the all-zero word) from $V''$; denote $F = \{ x \in V'\times V'' : \tilde h(x)=o \}$.
\begin{lemma}\label{l:face}
Assume that $G'$ is a regular graph of degree $d$.
Then the mapping $\tilde h$ defined by {\rm (\ref{eq:h})} is a perfect coloring of $G' \times G''$ with the parameter matrix $A'' + dI$ where
$A''$ is the adjacency matrix of $G''$ and $I$ is the identity matrix.
\end{lemma}

Arguing as in the previous subsection,
for any other perfect structure $f$ over $G' \times G''$ with parameters $S$
we have
$ (A'' + dI) (h\trn f) = (h\trn f) S, $
or, equivalently,
$$ A'' (h\trn f) = (h\trn f) (S-dI), $$
and, taking into account the obvious analog of Lemma~\ref{l:dist} for $F$, we derive the following:

\begin{theorem}\label{th:cartesian}
Let $f$ be an $S$-perfect structure over the cartesian product $G'\times G''$ of
a regular graph $G'=(V',E')$ and a distance-regular graph $G''=(V'',E'')$.
Let $F$ be some instance of $G'$ in $G'\times G''$, that is, the subgraph of $G'\times G''$
generated by the vertex subset $V' \times \{o\}$ for some $o\in V''$.
Let $f_0$ be the sum of $f$ over $F$.
Then the weight distribution of $f$ with respect to $F$ is
$$ (f_0=f_0\Pi_0(S-dI),\, f_0 \Pi_1(S-dI),\, \ldots ,\, f_0 \Pi_{\mathrm{diam}(G'')}(S-dI))\trn,$$
where $d$ is the degree of $G'$ and $\Pi_0$, $\Pi_1$, \ldots are the $P$-polynomials of $G''$.
\end{theorem}

A known example of the graph cartesian product is the $q$-ary $n$-cube $H^{n}_q$, which is
the cartesian product of $n$ copies of the full graph $K_q$.
For any integer $m$ from $0$ to $n$
we have $H^{n}_q = H^{m}_q \times H^{n-m}_q$, and Theorem~\ref{th:cartesian} means the following.

\begin{corollary}\label{cor:face}
Let $f$ be an $S$-perfect structure over the $q$-ary $(m+k)$-cube $H^{m+k}_q$.
Let $F$ be a subcube of $H^{m+k}_q$  of dimension $m$ (i.e., isomorphic to $H^{m}_q$),
and let $f_0$ be the sum of $f$ over $F$.
Then the weight distribution of $f$ with respect to $F$ is
$$ (f_0=f_0\Pi_0(S'),\, f_0 \Pi_1(S'),\, \ldots ,\, f_0 \Pi_{k}(S'))\trn,$$
where $S'=S-(q-1)mI$,
$\Pi_w(\cdot) = P_w(P_1^{-1}(\cdot))$,
and $P_w(x) = P_w(x;k,q)$ are the Krawtchouk polynomial {\rm (\ref{eq:krawtchuk})}.
The corresponding weight enumerator
equals $f_0 Z_{S';k,q}(z)$ where $Z_{\ldots}$ is defined in {\rm(\ref{eq:we0})}.
\end{corollary}

\subsection{A subcube of smaller size}
Here, we consider distributions with respect to another completely regular set in $H^{n}_q$
with large covering radius, the subcube $H^{n}_p$ of the same dimension and smaller order $p<q$.
Note that if $p$ divides $q$, then the theorem can be proved
using the approach of the previous two subsections; in the general case, we will calculate the parameters
of the distance coloring and use the known recursive formulas
\begin{eqnarray}
\textstyle (i+1)P_{i+1}(x;n,\frac qp)
&=&\textstyle((n-i)(\frac qp-1)+i-\frac qp x)P_i(x;n,\frac qp) \nonumber\\ &-& \textstyle(\frac qp-1)(n-i+1) P_{i-1}(x;n,\frac qp) \label{eq:Kqp}
\end{eqnarray}
for the Krawtchouk polynomials (see, e.g., \cite[\S\,5.7]{MWS}).
\begin{theorem}\label{th:pcube}
Let $f$ be an $S$-perfect structure over the $q$-ary $n$-cube $H^{n}_q$.
Let $p<q$, and let $f_0$ be the sum of $f$ over the
vertex of the $p$-ary subcube $H^{n}_p \subset H^{n}_q$.
Then the weight distribution of $f$ with respect to $V(H^{n}_p)$ is
$$ (f_0=f_0\Pi_0(S'),\, f_0 \Pi_1(S'),\, \ldots ,\, f_0 \Pi_{n}(S'))\trn,$$
where
$S'=(S-(p-1)nI)/p$, \ \ $\Pi_i(\cdot) = P_i(P_1^{-1}(\cdot))$, \ \
$P_i(x) = P_i(x;n,\frac qp)$ \ {\rm (\ref{eq:krawtchuk})}.
The corresponding weight enumerator
equals $f_0 Z_{\frac 1m S;k,\frac qp}(z)$ where $Z_{\ldots}$ is defined in {\rm(\ref{eq:we0})}.
\end{theorem}
\begin{proof}
The vertices of $H^{n}_q$ are the $n$-words over the alphabet $\{0,1,\ldots,q-1\}$,
while the vertices of $H^{n}_p$ are the $n$-words over the subalphabet $\{0,1,\ldots,p-1\}$.
Let us consider the distance coloring $g$ with respect to $V(H^{n}_p)$.
A word $v$ is at distance $w$ from $V(H^{n}_p)$, that is, $g(v)=e_w$,
if and only if it has $n-w$ symbols from
$\{0,1,\ldots,p-1\}$ and the other $w$ symbols from $\{p,\ldots,q-1\}$.
Every such word has $R_{w,w-1}=pw$ neighbors of color $e_{w-1}$,
$R_{w,w}=(p-1)(n-w)+(q-p-1)w$ neighbors of color $e_{w}$,
and $R_{w,w+1}=(q-p)(n-w)$ of color $e_{w+1}$.
Defining the other elements of the matrix $R$ as zeroes,
we obtain the parameter matrix
of the perfect coloring $g$. Now, we consider the transposed matrix $B=R\trn$.
Its nonzero elements are $B_{i,i-1}=(q-p)(n-i+1)$ (i.e., $R_{w,w+1}$ with $w=i-1$),
$B_{i,i}=(p-1)(n-i)+(q-p-1)i$, and $B_{i,i+1}=p(i+1)$ (i.e., $R_{w,w-1}$ with $w=i+1$).

Now consider an arbitrary $S$-perfect structure $f$ over $H^{n}_q$ and its distribution $h=g\trn f$
with respect to $g$.
By Theorems~\ref{th:distrib} and~\ref{th:main}, the rows $h_i$ of $h$ satisfy (\ref{eq:req3}),
that is, in our case,
$$
h_{i+1} = (h_{i} S - ((p-1)(n-i)+(q-p-1)i) h_{i} - (q-p)(n-i+1) h_{i-1})/p(i+1).
$$
Replacing $S$ by $pS'+(p-1)nI$ (i.e., $S'=(S-(p-1)nI)/p$), we get
$$
\textstyle(i+1)h_{i+1} = h_{i} S' - (\frac qp-2) i h_{i} - (\frac qp-1)(n-i+1) h_{i-1}.
$$
Now replace $\textstyle S'=P_1(S_P)=(\frac qp-1)n-\frac qp S_P$:
$$
\textstyle (i+1)h_{i+1} = -\frac qp h_{i}  S_P +((\frac qp-1)n - (\frac qp-2) i) h_{i} - (\frac qp-1)(n-i+1) h_{i-1}.
$$
Since this recursive relation coincides with (\ref{eq:Kqp}) and the first
elements $h_0$, $h_1$ of the sequence also coincide with $h_0 P_0(S_P)$, $h_0 P_1(S_P)$,
we see that $h_i = h_0 P_i(S_P)=h_0 P_i(P_1^{-1}(S'))$ for every $i$.
\qed\end{proof}
\section{Local distributions in the cartesian product of graphs.}\label{s:local}
Let us consider two graphs $G'$ and $G''$
and select one vertex in every graph,
say, $o'\in V(G')$ and $o''\in V(G'')$.
Assume that the distance colorings $g'$ and $g''$ with respect to $o'$ and $o''$, respectively, are perfect.
Consider some perfect coloring $f$ of the cartesian product $G=G'\times G''$.
It generates some colorings (not necessarily perfect) $f'$ and $f''$ of the subgraphs
$G' \times o''$ and $o' \times G''$, which are isomorphic to $G'$ and $G''$, respectively.
The distributions of $f'$ and $f''$ with respect to the vertex $(o',o'')$
(in the graphs $G' \times o''$ and $o' \times G''$, respectively)
will be called \emph{local distributions} (local spectra \cite{Vas04:local,Vas09:inter}) of $f$.
It turns out, one of the local distributions (say, of $f''$)
uniquely defines the other (of $f'$) \cite{Avg:lect}.
This fact was first proved in \cite{Vas04:local} for $1$-perfect codes in binary $n$-cubes;
an explicit formula was derived, see also \cite[Th.~3]{AvgVas:2006:testing}.
Our goal is to derive a matrix formula that connects the local distributions
(more tightly, the formula for the distribution of $f$ with
respect to the perfect coloring $g'\otimes g''$;
this distribution includes the local distributions).
We start from the general case, when $G'$ and $G''$ are arbitrary graphs, $g'$ and $g''$ are
arbitrary perfect colorings (or, even more generally, perfect structures).
As in the previous sections, in partial cases, the formula will have explicit solutions.

\begin{figure}
\noindent\mbox{}\hfill
\includegraphics[scale=1.6]{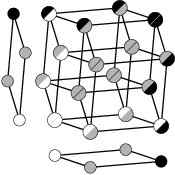}
\hfill
\small\raisebox{6em}{$\left(
 \begin{array}{cc@{\ \,}c@{\ \,}c@{\ \ \ }c@{\ \,}c@{\ \,}c@{\ \ \ }c@{\ \,}c@{\ \,}c}
 &
 \includegraphics{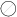} &
 \includegraphics{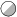} &
 \includegraphics{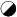} &
 \includegraphics{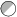} &
 \includegraphics{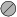} &
 \includegraphics{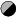} &
 \includegraphics{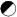} &
 \includegraphics{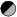} &
 \includegraphics{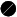}
 \\
 \includegraphics{mmm100} & 0 & 2 & 0 & 2 & 0 & 0 & 0 & 0 & 0 \\[-0.7mm]
 \includegraphics{mmm101} & 1 & 0 & 1 & 0 & 2 & 0 & 0 & 0 & 0 \\[-0.7mm]
 \includegraphics{mmm102} & 0 & 2 & 0 & 0 & 0 & 2 & 0 & 0 & 0 \\[1.1mm]
 \includegraphics{mmm110} & 1 & 0 & 0 & 0 & 2 & 0 & 1 & 0 & 0 \\[-0.7mm]
 \includegraphics{mmm111} & 0 & 1 & 0 & 1 & 0 & 1 & 0 & 1 & 0 \\[-0.7mm]
 \includegraphics{mmm112} & 0 & 0 & 1 & 0 & 2 & 0 & 0 & 0 & 1 \\[1.1mm]
 \includegraphics{mmm120} & 0 & 0 & 0 & 2 & 0 & 0 & 0 & 2 & 0 \\[-0.7mm]
 \includegraphics{mmm121} & 0 & 0 & 0 & 0 & 2 & 0 & 1 & 0 & 1 \\[-0.7mm]
 \includegraphics{mmm122} & 0 & 0 & 0 & 0 & 0 & 2 & 0 & 2 & 0
 \end{array}
 \right)$}
 \hfill\mbox{}
\caption{\label{fig:H22}Cartesian product $H_2^2\times H_2^2$ and tensor product of distance colorings}
\end{figure}

We first consider the representation of the cartesian product of graphs by its
adjacency matrix.
The \emph{tensor product} $A' \otimes A''$ of $n' \times m'$ and
$n'' \times m''$ matrices $A'=(a'_{ij})$ and $A''=(a''_{ij})$ is defined as the
$n'n'' \times m'm''$ matrix $A=(a_{i'i''j'j''})$ whose
rows are indexed by two numbers
$i'=0,\ldots,n'-1$ and $i''=0,\ldots,n''-1$,
columns are indexed by two numbers
$j'=0,\ldots,m'-1$ and $j''=0,\ldots,m''-1$,
and the elements $a_{i'i''j'j''}$ are equal to $a'_{i'j'}a''_{i''j''}$.
The following well-known fact is straightforward from the definitions.
\begin{lemma}
  The adjacency matrices $A'$, $A''$, and $A$ of graphs $G'$, $G''$, and their cartesian product $G = G'\times G''$, respectively,
  are related by $A=A'\otimes I + I \otimes A''$.
\end{lemma}
From simple counting arguments, we have the following:
\begin{lemma}\label{l:R'R''}
Let $g'$ and $g''$ be respectively $R'$- and $R''$- perfect structures over graphs $G'$ and $G''$.
Then $g'\otimes g''$ is an $(R'\otimes I + I\otimes R'')$-perfect structure over the cartesian product $G'\times G''$.
Briefly,
\begin{multline*}
(A'g'=g'R') \ \&\  (A''g''=g''R'') \\
 \Rightarrow ((A'\otimes I + I\otimes A'')(g'\otimes g'')=(g'\otimes g'')(R'\otimes I + I\otimes R'')).
\end{multline*}
\end{lemma}
\begin{proof}
The implication follows immediately from the straightforward property
$ (X\otimes Y)(Z\otimes V) = (XZ)\otimes (YV) $ of the tensor product.
\qed\end{proof}
\begin{remark}
  Assume that $g'$ and $g''$ are perfect colorings, that is, every row contains $1$ in one position
  and $0$s in the others. Then $g=g' \otimes g''$ satisfies the same property.
  Indeed, $g_{i'i''j'j''}=1$ if and only if $g'_{i'j'}=1$ and $g''_{i''j''}=1$.
\end{remark}
\begin{example}
  Let $G'=G''$ be the $3$-ary $2$-cube. Then $G=G' \times G''$ is the $3$-ary $4$-cube.
  Let $g'$ be the distance coloring with respect to some point in $G'$; $g''$, in $G''$. Then
  the parameter matrix of the perfect coloring $g' \otimes g''$ is the following:
 $$ \left(
 \begin{array}{ccc}
 0 & 4 & 0 \\
 1 & 1 & 2 \\
 0 & 2 & 2
 \end{array}
 \right)
 \otimes I + I \otimes
 \left(
 \begin{array}{ccc}
 0 & 4 & 0 \\
 1 & 1 & 2 \\
 0 & 2 & 2 \end{array}
 \right)
 =
 \left(
 \begin{array}{c@{\ \,}c@{\ \,}c@{\ \ \ }c@{\ \,}c@{\ \,}c@{\ \ \ }c@{\ \,}c@{\ \,}c}
 0 & 4 & 0 & 4 & 0 & 0 & 0 & 0 & 0 \\[-0.7mm]
 1 & 1 & 2 & 0 & 4 & 0 & 0 & 0 & 0 \\[-0.7mm]
 0 & 2 & 2 & 0 & 0 & 4 & 0 & 0 & 0 \\[1.1mm]
 1 & 0 & 0 & 1 & 4 & 0 & 2 & 0 & 0 \\[-0.7mm]
 0 & 1 & 0 & 1 & 2 & 2 & 0 & 2 & 0 \\[-0.7mm]
 0 & 0 & 1 & 0 & 2 & 3 & 0 & 0 & 2 \\[1.1mm]
 0 & 0 & 0 & 2 & 0 & 0 & 2 & 4 & 0 \\[-0.7mm]
 0 & 0 & 0 & 0 & 2 & 0 & 1 & 3 & 2 \\[-0.7mm]
 0 & 0 & 0 & 0 & 0 & 2 & 0 & 2 & 4
 \end{array}
 \right)
 $$
For a smaller example of the cartesian product of two binary $2$-cubes see Fig.~\ref{fig:H22}.
\end{example}
We consider the distribution of some $S$-perfect structure (perfect coloring) $f$ over $G=G'\times G''$
with respect to the perfect structure $g = g'\otimes g''$, that is, the matrix
$h=(g'\otimes g'')\trn f = ({g'}\trn\otimes {g''}\trn) f$.
The rows $h_{i'i''}$ of $h$ are indexed by two indices,
corresponding to the colors of $g'$ and $g''$ respectively;
the columns are indexed by the colors of $f$.
By Lemma~\ref{l:R'R''} and Theorem~\ref{th:distrib}, $h$ satisfies 
\begin{equation}\label{eq:RIIR}
({R'}\trn \otimes I + I \otimes {R''}\trn)h = h S
\end{equation}
Our goal is, provided $g'$ is a distance coloring with respect to some set $C$
(e.g., $C=\{c\}$),
to reconstruct $h$ from knowledge of only rows of type $h_{0i''}$,
that is, from knowledge of the distribution
of $f$ with respect to the restriction of $g$ to $C\times G''$ 
(if $C$ consists of one vertex, then this restriction
is isomorphic to the coloring $g''$ of $G''$).
To do this, we rearrange the elements of the matrix $h$ in such a way that
all known elements are in the first row of the new matrix, say $h^*$.
The element $h_{i'i''j}$ of $h$ coincides with the corresponding element of
$h^*$, but in $h^*$, the first index is the row number,
while the second and the third index the columns.
So, if $h$ is a ${\mu}'{\mu}''\times {\nu}$ matrix,
then $h^*$ is a ${\mu}'\times {\mu}''{\nu}$ matrix.
With $h^*$, the equation (\ref{eq:RIIR}) can be rewritten as follows:
$${R'}\trn h^* + h^*({R''}\otimes I) = h^* (I\otimes S)$$
or
$${R'}\trn h^* = h^*(I\otimes S-{R''}\otimes I)$$
So, we have proved the following:
\begin{theorem}
Let $g'$, $g''$, and $f$ be $R'$-, $R''$-, and $S$- perfect colorings of
graphs $G'$, $G''$, and $G = G' \times G''$ respectively.
Let $h = (g' \otimes g'')\trn f$ be the distribution of $f$ with respect to
the perfect coloring $g' \otimes g''$ of $G$.
Then $h^*$ is an $( I\otimes S-{R''}\otimes I)$-perfect structure over ${R'}\trn$.
\end{theorem}
\begin{corollary}
  If $g'$ is a distance coloring with respect to a vertex in a distance-regular graph $G'$ with
  P-polynomials $\Pi_0$, $\Pi_1$, \ldots, then the rows $h^*_i$ of $h^*$ can be calculated as
    \begin{equation}\label{eq:loc}
   h^*_i = h^*_0 \Pi_i(I\otimes S-{R''}\otimes I).
\end{equation}
\end{corollary}
\begin{remark}
  If $g'$ ($g''$) is a distance colorings, then the submatrix $(h_{0i''j})$
  (respectively, $(h_{i'0j})$) of $h=(h_{i'i''j})$ is a local distribution of $f$,
  see the introduction of the section.
\end{remark}
\begin{remark}
  If $g'$ is a \emph{trivial} perfect coloring (each vertex is colored into its own color),
  then $R'$ coincides with the adjacency matrix of $G'$,
  and $h^*$ is a perfect structure over $G'$.
\end{remark}
\begin{example}
  Let $g'$ and $g''$ be the distance colorings with respect to the all-zero word in the Hamming
  graphs $H_q^{n-3}$ and $H_q^3$ respectively. The perfect coloring $g''$ has the parameters
  $$
  R''=\left(\begin{array}{cccc}0&3q{-}3&0&0\\1&q{-}2&2q{-}2&0\\0&2&2q{-}4&q{-}1\\0&0&3&3q{-}6\end{array}\right).
  $$
  Let $f$ be an $S$-perfect coloring with $S=\left(\begin{array}{cc}0&(q{-}1)n\\1&(q{-}1)n{-}1\end{array}\right)$,
  that is, the first color corresponds to a $1$-perfect code $C$, see Example~\ref{ex:1perf}.
  Then $
  I\otimes S - R'' \otimes I = $ {\small$$
   \left(
 \begin{array}{c@{\ \,}c@{\ |\ }c@{\ \,}c@{\ |\ }c@{\ \,}c@{\ |\ }c@{\ \,}c}
 0       & (q{-}1)n& -3q{+}3 & 0       & 0       & 0       & 0       & 0       \\
 1    &(q{-}1)n{-}1& 0       & -3q{+}3 & 0       & 0       & 0       & 0       \\
\hline
 -1      & 0       & -q{+}2  & (q{-}1)n& -2q{+}2 & 0       & 0       & 0       \\
 0       & -1      & 1 & (q{-}1)(n{-}1)& 0       & -2q{+}2 & 0       & 0       \\
\hline
 0       & 0       & -2      & 0       & -2q{+}4 & (q{-}1)n& -q{+}1  & 0       \\
 0       & 0       & 0       & -2      & 1&(q{-}1)(n{-}2){+}1& 0       & -q{+}1  \\
\hline
 0       & 0       & 0       & 0       & -3      & 0       & -3q{+}6 & (q{-}1)n\\
 0       & 0       & 0       & 0       & 0       & -3      & 1       & (q{-}1)(n{-}3){+}2
 \end{array}
 \right)
 $$}
\end{example}
Now suppose that the code $C$ contains the all-zero word.
Since the minimal distance between codewords is $3$, the only possibilities for
$h^*_0$ are
$$(1,0;\ 0,3q{-}3;\ 0,3(q{-}1)^2;\ t, (q{-}1)^3{-}t)\quad t\in \{0,1,\ldots, q-1\},$$
that is, there are from $0$ to $q-1$ nonzero codewords in $H_q^3\times 0^{n-3}$,
and all of them are at distance $3$ from the all-zero word. Further, if $n=q+1$, then $t=q-1$
(this comes from the fact that, by numerical reasons, every subgraph isomorphic to $H_q^2$ contains exactly one code vertex; as follows, $H_q^3\times 0^{n-3}$ contains exactly $q$ code vertices).
Substituting $q=10$ and $n=11$ and calculating $h^*$ by formulas (\ref{eq:loc}), (\ref{eq:krawtchuk}), we get
{
\footnotesize
$$
   \left(
 \begin{array}{c@{\ \,}c@{\ |\ \scriptstyle}c@{\ \,}c@{\ |\ }c@{\ \,}c@{\ |\ }c@{\ \,}c}
1 & \scriptstyle 0 &    0 & \scriptstyle 27 & 0 & \scriptstyle 243 & 9 & \scriptstyle 720   \\
0 & \scriptstyle 72 & 0 & \scriptstyle 1944 & 216 & \scriptstyle 17280 & 504 & \scriptstyle 51984    \\
0 & \scriptstyle 2268 & 756 & \scriptstyle 60480 & 5292 & \scriptstyle 545832 & 16632 & \scriptstyle 1636740    \\
504 & \scriptstyle 40320 & 10584 & \scriptstyle 1091664 & 99792 & \scriptstyle 9820440 & 297360 & \scriptstyle 29463336 \\
4410 & \scriptstyle 454860 & 124740 & \scriptstyle 12275550 & 1115100 & \scriptstyle 110487510 & 3348450 & \scriptstyle 331459380   \\
33264 & \scriptstyle 3273480 & 892080 & \scriptstyle 88390008 & 8036280 & \scriptstyle 795502512 & 24105816 & \scriptstyle 2386510560   \\
148680 & \scriptstyle 14731668 & 4018140 & \scriptstyle 397751256 & 36158724 & \scriptstyle 3579765840 & 108477936 & \scriptstyle 10739295756   \\
382680 & \scriptstyle 37881072 & 10331064 & \scriptstyle 1022790240 & 92981088 & \scriptstyle 9205110648 & 278942688 & \scriptstyle 27615332520 \\
430461 & \scriptstyle 42616260 & 11622636 & \scriptstyle 1150638831 & 104603508 & \scriptstyle 10355749695 & 313810605 & \scriptstyle 31067249004
 \end{array}\right)
$$
}
So, if a $1$-perfect code exists in $H_{10}^{11}$
(this is an open question) and contains the all-zero word,
its distribution with respect to $g' \otimes g''$
is given by the matrix above.

\section{Conclusions}
We have derived quite general matrix formulas for different weight distributions of
perfect structures (perfect colorings, completely regular codes).
One of the interesting open problems in this topic is to obtain formulas
for the interweight distribution of a perfect coloring in a binary $n$-cube $H_2^n$.
By the \emph{interweight distribution} (interweight spectrum \cite{Vas09:inter}) of an $S$-perfect coloring,
with respect to a point $a$,
we mean the set of values
\begin{eqnarray*}
T^{st}_{ijd}(a) = |\{(x, y) & : & x,y\in V(H_2^n),\
f(x) = s,\ f(y) = t, \\
&&d_{H_2^n}(x, a) = i,\ d_{H_2^n}(y, a) = j,\ d_{H_2^n}(x, y) = d\}|.
\end{eqnarray*}
As was found in \cite{Vas09:inter} (using the local distributions
defined in Section~\ref{s:local}), the interweight distribution
depends only on the parameter matrix $A$ and the color of the initial point $a$, and does not depend on the choice of the perfect coloring and the initial point.
This cannot be generalized to an arbitrary distance-regular graph.
For example, the two sets
\begin{eqnarray*}
&&\{000000,110000,111100,111111,001111,000011\},\\
&&\{000000,110000,011000,111111,001111,100111\}
\end{eqnarray*}
 in the halved $6$-cube $H_+^6$ both have $S$-perfect distance colorings with the same $S$; but
$T^{00}_{111}(000000)=0$ in the first case and
$T^{00}_{111}(000000)=1$ in the second.
Another example is the union of the three $1$-perfect codes
\begin{eqnarray*}
&&\{0000,0111,0222,1012,1120,1201,2021,2102,2210\},\\
&&\{1000,1111,1222,2012,2120,2201,0021,0102,0210\},\\
&&\{2000,0011,0202,1022,0120,1101,2221,2112,1210\}
\end{eqnarray*}
in $H_3^4$,
whose distance coloring is perfect with
$1=T^{00}_{111}(0000)\neq T^{00}_{111}(0111)=0$.

For binary $n$-cubes, it would be interesting to derive an invariant of perfect structures that generalizes the interweight distribution of perfect colorings.


\begin{thebibliography}{10}
\providecommand{\url}[1]{{#1}}
\providecommand{\urlprefix}{URL }
\expandafter\ifx\csname urlstyle\endcsname\relax
  \providecommand{\doi}[1]{DOI~\discretionary{}{}{}#1}\else
  \providecommand{\doi}{DOI~\discretionary{}{}{}\begingroup
  \urlstyle{rm}\Url}\fi

\bibitem{Avg:PhD}
Avgustinovich, S.V.: Metrical and combinatorial properties of perfect codes and
  colorings.
\newblock {PhD} thesis, Sobolev Institute of Mathematics, Novosibirsk, Russia
  (2000).
\newblock In Russian

\bibitem{Avg:lect}
Avgustinovich, S.V.: Perfect structures (2007).
\newblock Lectures. Korea, POSTECH. Unpublished

\bibitem{AvgMog:J63J73}
Avgustinovich, S.V., Mogilnykh, I.Y.: Perfect $2$-colorings of {J}ohnson graphs
  ${J}(6,3)$ and ${J}(7,3)$.
\newblock In: \'{A}ngela Barbero (ed.) Coding Theory and Applications (Second
  International Castle Meeting, {ICMCTA} 2008, Castillo de la Mota, Medina del
  Campo, Spain, September 15-19, 2008. Proceedings), \emph{Lect. Notes Comput.
  Sci.}, vol. 5228, pp. 11--19. Springer-Verlag, Berlin Heidelberg (2008).
\newblock \DOI{10.1007/978-3-540-87448-5\_2}

\bibitem{AvgVas:2006:testing}
Avgustinovich, S.V., Vasil'eva, A.Y.: Testing sets for 1-perfect code.
\newblock In: R.~Ahlswede, L.~B{\"a}umer, N.~Cai, H.~Aydinian, V.~Blinovsky,
  C.~Deppe, H.~Mashurian (eds.) General Theory of Information Transfer and
  Combinatorics, \emph{Lect. Notes Comput. Sci.}, vol. 4123, pp. 938--940.
  Springer-Verlag, Berlin Heidelberg (2006).
\newblock \DOI{10.1007/11889342\_59}

\bibitem{BZZ:1974:UPC}
Bassalygo, L.A., Zaitsev, G.V., Zinoviev, V.A.: On uniformly packed codes.
\newblock \href{http://www.springerlink.com/content/0032-9460}{Probl. Inf.
  Transm.} \textbf{10}(1), 6--9 (1974).
\newblock Translated from \em Probl. Peredachi Inf.\em, 10(1): 9-14, 1974

\bibitem{Brouwer}
Brouwer, A.E., Cohen, A.M., Neumaier, A.: Distance-Regular Graphs.
\newblock Springer-Verlag, Berlin (1989)

\bibitem{CDS}
Cvetkovi\'c, D.M., Doob, M., Sachs, H.: Algebraic Graph Theory.
\newblock Press Inc., New York (1980)

\bibitem{Delsarte:1973}
Delsarte, P.: An algebraic approach to association schemes of coding theory.
\newblock Philips Research Reports, Supplement \textbf{10} (1973)

\bibitem{FDF:PerfCol}
Fon-Der-Flaass, D.G.: Perfect $2$-colorings of a hypercube.
\newblock \href{http://www.springerlink.com/content/0037-4466}{Sib. Math. J.}
  \textbf{48}(4), 740--745 (2007).
\newblock \DOI{10.1007/s11202-007-0075-4} translated from
  \href{http://math.nsc.ru/smz/}{Sib. Mat. Zh.} 48(4) (2007), 923-930

\bibitem{FDF:12cube.en}
Fon-Der-Flaass, D.G.: Perfect colorings of the $12$-cube that attain the bound
  on correlation immunity.
\newblock Sib. Ehlektron. Mat. Izv. \textbf{4}, 292--295 (2007).
\newblock In Russian. Online: \url{http://semr.math.nsc.ru/v4/p292-295.pdf}

\bibitem{Godsil93}
Godsil, C.D.: Algebraic Combinatorics.
\newblock Chapman and Hall, New York (1993)

\bibitem{Hed:2008:survey}
Heden, O.: A survey of perfect codes.
\newblock \href{http://aimsciences.org/journals/amc/}{Adv. Math. Commun.}
  \textbf{2}(2), 223--247 (2008). \newblock \DOI{10.3934/amc.2008.2.223}

\bibitem{Kro:24}
Krotov, D.S.: On perfect colorings of the halved $24$-cube.
\newblock \href{http://www.math.nsc.ru/publishing/DAOR/daor.html}{Diskretn.
  Anal. Issled. Oper.} \textbf{15}(5), 35--46 (2008).
\newblock In Russian; translated at \url{http://arxiv.org/abs/0803.0068}

\bibitem{Kro:2m-3}
Krotov, D.S.: On the binary codes with parameters of doubly-shortened
  $1$-perfect codes.
\newblock \href{http://www.springerlink.com/content/100256/}{Des. Codes
  Cryptography} \textbf{57}(2), 181--194 (2010).
\newblock \DOI{10.1007/s10623-009-9360-5}

\bibitem{Kro:2010ACCT:n-3}
Krotov, D.S.: On the binary codes with parameters of triply-shortened
  $1$-perfect codes.
\newblock Submitted. arXiv: \href{http://arXiv.org/abs/1104.0005}{1104.0005}

\bibitem{Lloyd}
Lloyd, S.P.: Binary block coding.
\newblock Bell Syst. Tech. J. \textbf{36}(2), 517--535 (1957)

\bibitem{MWS}
MacWilliams, F.J., Sloane, N.J.A.: The Theory of Error-Correcting Codes.
\newblock Amsterdam, Netherlands: North Holland (1977)

\bibitem{Martin:PhD}
Martin, W.J.: Completely regular subsets.
\newblock {PhD} thesis, University of Waterloo, Waterloo, Ontario, Canada
  (1992).
\newblock Online: \href{http://users.wpi.edu/~martin/RESEARCH/THESIS/}%
  {http://users.wpi.edu/{\char126}martin/RESEARCH/THESIS/}

\bibitem{Mar:CRD}
Martin, W.J.: Completely regular designs.
\newblock \href{http://www3.interscience.wiley.com/cgi-bin/jhome/38682}{J.
  Comb. Des.} \textbf{6}(4), 261--273 (1998).
\newblock
  \DOIURL{10.1002/(SICI)1520-6610(1998)6:4$<$261::AID-JCD4$>$3.0.CO;2-D}{10.10%
02/(SICI)1520-6610(1998)6:4<261::AID-JCD4>3.0.CO;2-D}

\bibitem{SZZ:1971:UPC}
Semakov, N.V., Zinoviev, V.A., Zaitsev, G.V.: Uniformly packed codes.
\newblock \href{http://www.springerlink.com/content/0032-9460}{Probl. Inf.
  Transm.} \textbf{7}(1), 30--39 (1971).
\newblock Translated from \em Probl. Peredachi Inf.\em, 7(1): 38-50, 1971

\bibitem{ShSl}
Shapiro, H.S., Slotnick, D.L.: On the mathematical theory of error correcting
  codes.
\newblock IBM J. Res. Develop. \textbf{3}(1), 25--34 (1959)

\bibitem{Sol:2008:survey}
Solov'eva, F.I.: On perfect binary codes.
\newblock \href{http://www.sciencedirect.com/science/journal/0166218X}{Discrete
  Appl. Math.} \textbf{156}(9), 1488--1498 (2008).
\newblock \DOI{10.1016/j.dam.2005.10.023}

\bibitem{Vas:cr2}
Vasil'eva, A.: Linear binary completely regular codes with distance $2$.
\newblock In: Proc. 2008 {IEEE} {R}egion 8 International Conference on
  Computational Technologies in Electrical and Electronics Engineering
  ``SIBIRCON 2008'', pp. 12--15. Novosibirsk, Russia (2008).
\newblock \DOI{10.1109/SIBIRCON.2008.4602615}

\bibitem{Vas04:local}
Vasil'eva, A.Y.: Local spectra of perfect binary codes.
\newblock \href{http://www.sciencedirect.com/science/journal/0166218X}{Discrete
  Appl. Math.} \textbf{135}(1-3), 301--307 (2004).
\newblock \DOI{10.1016/S0166-218X(02)00313-X}, translated from
  \href{http://www.math.nsc.ru/publishing/DAOR/daor.html}{Diskretn. Anal.
  Issled. Oper.}, Ser.~1, 6(1):3-11, 1999

\bibitem{Vas09:inter}
Vasil'eva, A.Y.: Local and interweight spectra of completely regular codes and
  of perfect colorings.
\newblock \href{http://www.springerlink.com/content/0032-9460}{Probl. Inf.
  Transm.} \textbf{45}(2), 151--157 (2009).
\newblock \DOI{10.1134/S0032946009020069}, translated from Probl. Peredachi
  Inf., 45(2):84-90, 2009

\bibitem{ZinHel04}
Zinoviev, V.A., Helleseth, T.: On weight distributions of shifts of
  {G}oethals-like codes.
\newblock \href{http://www.springerlink.com/content/0032-9460}{Probl. Inf.
  Transm.} \textbf{40}(2), 118--134 (2004).
\newblock \DOI{10.1023/B:PRIT.0000043926.60991.7e} translated from Probl.
  Peredachi Inf. 40(2) (2004), 19-36

\bibitem{RifZin07}
Zinoviev, V.A., Rif\`{a}, J.: On new completely regular $q$-ary codes.
\newblock \href{http://www.springerlink.com/content/0032-9460}{Probl. Inf.
  Transm.} \textbf{43}(2), 97--112 (2007).
\newblock \DOI{10.1134/S0032946007020032} translated from Probl. Peredachi Inf.
  43(2) (2007), 34-51

\end{thebibliography}

\providecommand\href[2]{#2} \providecommand\url[1]{\href{#1}{#1}}
  \def\DOI#1{{DOI}: \href{http://dx.doi.org/#1}{#1}}
	\def\DOIURL#1#2{ {DOI}:\href{http://dx.doi.org/#2}{#1}}

\end{document}